\documentclass[11pt]{amsart}

\usepackage{amsmath,amssymb,amsthm,amsrefs,a4wide}
\usepackage{latexsym}
\usepackage{indentfirst}
\usepackage{graphicx}
\usepackage{placeins}
\usepackage{booktabs}
\usepackage{algorithm}
\usepackage{algorithmic}
\usepackage{multirow}
\usepackage{color}
\usepackage{esint}

\theoremstyle{plain}

\theoremstyle{definition}

\newtheorem{example}{Example}

\theoremstyle{remark}
\newtheorem{remark}{Remark}

\newcommand{\diag}{\mathrm{diag}}

\newcommand{\bbm}{\begin{bmatrix}}
\newcommand{\ebm}{\end{bmatrix}}
\newcommand{\R}{\mathbb{R}}
\newcommand{\C}{\mathbb{C}}
\newcommand{\D}{\mathbb{D}}

\newcommand{\nx}{{n_x}}
\newcommand{\ns}{{n_s}}
\newcommand{\na}{{n_a}}
\renewcommand{\L}{\ell}

\newcommand{\ut}{\tilde{u}}
\newcommand{\xt}{\tilde{x}}
\newcommand{\wt}{\tilde{w}}

\newcommand{\Zt}{\tilde{Z}}

\newcommand{\Gh}{\hat{G}}

\newcommand{\bg}{\mathbf{g}}

\newcommand{\bu}{\mathbf{u}}

\newcommand{\but}{\tilde{\bu}}

\newcommand{\bgh}{\hat{\bg}}

\begin{document}

\title[Multidimensional Eigenmatrix]{Multidimensional unstructured sparse recovery via eigenmatrix}

\author[]{Lexing Ying} \address[Lexing Ying]{Department of Mathematics, Stanford University,
  Stanford, CA 94305} \email{lexing@stanford.edu}

\thanks{This work is partially supported by NSF grants DMS-2011699 and DMS-2208163.}

\keywords{Sparse recovery, Prony's method, ESPRIT, eigenmatrix.}

\subjclass[2010]{30B40, 65R32.}

\begin{abstract}
  This note considers the multidimensional unstructured sparse recovery problems. Examples include Fourier inversion and sparse deconvolution. The eigenmatrix is a data-driven construction with desired approximate eigenvalues and eigenvectors proposed for the one-dimensional problems. This note extends the eigenmatrix approach to multidimensional problems. Numerical results are provided to demonstrate the performance of the proposed method.
\end{abstract}

\maketitle

\section{Introduction}\label{sec:intro}

This note considers the multidimensional unstructured sparse recovery problems of a general form. Let $X$ be the parameter space, typically a subset of $\R^d$ or $\C^d$, and $S$ be the sampling space. $G(s,x)$ is a kernel function defined for $s\in S$ and $x\in X$, and is assumed to be analytic in $x$. Suppose that
\begin{equation}
  f(x) = \sum_{k=1}^{\nx} w_k \delta(x-x_k)
  \label{eq:f}
\end{equation}
is the unknown multidimensional sparse signal, where $\{x_k\}_{1\le k \le n_x}$ are the spikes and $\{w_k\}_{1\le k \le n_x}$ are the weights. The observable of the problem is
\begin{equation}
  u(s) :=\int_X G(s,x) f(x) dx = \sum_{k=1}^{\nx} G(s,x_k) w_k
  \label{eq:u}
\end{equation}
for $s\in S$.

Let $\{s_j\}$ be a set of $\ns$ unstructured samples in $S$ and $u_j:=u(s_j)$ be the exact values. Suppose that we are only given the noisy observations $\ut_j:=u_j (1+\sigma Z_j)$, where $Z_j$ are independently identically distributed (i.i.d.) random variables with zero mean and unit variance, and $\sigma$ is the noise magnitude. The task is to recover the spikes $\{x_k\}$ and weights $\{w_k\}$ from $\{\ut_j\}$.

This note focuses on the case $X\subset \R^d$ since most applications fit into this real case. Two important examples are listed as follows.
\begin{itemize}
\item Sparse deconvolution. For example, $G(s,x) = \frac{1}{\|s-x\|^\alpha}$ or $e^{-\|s-x\|}$, $X$ is a bounded domain, and $\{s_j\}$ are samples outside $X$ in $\R^d$.
\item Fourier inversion. For example $G(s,x) = \exp(\pi i s\cdot x)$, $X=[-1,1]^d$ for example, and $\{s_j\}$ are samples in $\R^d$.
\end{itemize}

The primary challenges of the current setup come from three sources. First, the kernel $G(s,x)$ can be quite general. Second, the samples $\{s_j\}$ are unstructured, which excludes the existing algorithms that exploit the Cartesian structure. Third, the sample values $\{\ut_j\}$ are noisy, which raises stability issues when the recovery problem is quite ill-posed.

\subsection{Contribution}
The recent paper \cite{ying2023eigenmatrix} introduces the eigenmatrix for one-dimensional unstructured sparse recovery problems. This is a data-driven object that depends on $G(\cdot,\cdot)$, $X$, and the samples $\{s_j\}$. Constructed to have desired eigenvalues and eigenvectors, it turns the unstructured sparse recovery problem into an algebraic problem (such as a rootfinding or eigenvalue problem). The main features of the eigenmatrix are
\begin{itemize}
\item It assumes no special structure of the samples $\{s_j\}$.
\item It offers a rather unified approach to these sparse recovery problems.
\end{itemize}
This note extends this approach to the multidimensional cases.


\subsection{Related work} There has been a long list of works devoted to the sparse recovery problems mentioned above. We refer to the paper \cite{ying2023eigenmatrix} for the related work on the one-dimensional case. For the multidimensional Fourier inversion or superresolution, one well-studied approach is based on convex relaxation or $\ell_1$-minimization; we refer to \cite{poon2019multidimensional} and the references therein. There is a long list of works on the Cartesian samples based on Prony's method and the ESPRIT algorithm \cite{andersson2018esprit,andersson2010nonlinear,cuyt2011sparse,harmouch2018structured,kunis2016multivariate,peter2015prony,potts2013parameter,rouquette2001estimation,sacchini1993two,sahnoun2017multidimensional,sauer2017prony,vanpoucke1994efficient}. A closely related problem is the sparse moment problem. Among the references on this topic, \cite{fan2023efficient} is also closely related to our method for the 2D case.

The rest of the note is organized as follows. Section \ref{sec:review} reviews the eigenmatrix for one-dimensional problems in the setting of the ESPRIT algorithm. Section \ref{sec:method} describes the extension to multidimensional problems. Section \ref{sec:num} presents several numerical experiments. Section \ref{sec:disc} concludes with a discussion for future work.

\section{Eigenmatrix with ESPRIT}\label{sec:review}

This section provides a short review of the eigenmatrix approach. It was proposed in \cite{ying2023eigenmatrix} for both Prony's method \cite{prony1795essai} and the ESPRIT algorithm \cite{roy1989esprit}. The presentation here focuses on ESPRIT since it will be used for multidimensional problems.

{\bf Complex analytic case.} To simplify the discussion, assume that $X$ is the unit disc $\D\subset\C$. Define for each $x$ the vector $\bg(x):=[G(s_j,x)]_{1\le j\le \ns}$. The first step is to construct $M$ such that $M \bg(x) \approx x \bg(x)$ for $x\in \D$. Numerically, it is more robust to use the normalized vector $\bgh(x) = \bg(x)/\|\bg(x)\|$ since the norm of $\bg(x)$ can vary significantly depending on $x$. The condition then becomes
\[
M \bgh(x) \approx x \bgh(x), \quad x\in \D.
\]
We enforce this condition on a uniform grid $\{a_t\}_{1\le t \le \na}$ of size $\na$ on the boundary of the unit disk
\[
M \bgh(a_t) \approx a_t \bgh(a_t).
\]
Define the $\ns\times \na$ matrix $\Gh = [\bgh(a_t)]_{1\le t \le \na}$ with $\bgh(a_t)$ as columns and also the $\na\times \na$ diagonal matrix $\Lambda = \diag(a_t)$. The above condition can be written in a matrix form as
\[
M \Gh \approx \Gh \Lambda.
\]
$n_a$ is chosen so that the columns of $\Gh$ are numerically linearly independent (in practice, the condition number of $\Gh$ is bounded below $10^7$). When the columns of $\Gh$ are numerically linearly independent, we define the eigenmatrix as
\begin{equation}
  M := \Gh \Lambda \Gh^+, \label{eq:M}
\end{equation}
where the pseudoinverse $\Gh^+$ is computed by thresholding the singular values of $\Gh$. In practice, the thresholding value is chosen so that the norm of $M$ is bounded by a small constant.

{\bf Real analytic case.} To simplify the discussion, assume that $X$ is the interval $[-1,1]$. For each $x$, again $\bg(x) = [G(s_j,x)]_{1\le j\le \ns}$. The first step is to construct $M$ such that $M \bg(x) \approx x \bg(x)$ for $x\in [-1,1]$. Moving to the normalized vector $\bgh(x) = \bg(x)/\|\bg(x)\|$, we instead aim for
\[
M \bgh(x) \approx x \bgh(x), \quad x\in [-1,1].
\]
This condition is enforced on a Chebyshev grid $\{a_t\}_{1\le t \le \na}$ of size $\na$ on the interval $[-1,1]$:
\[
M \bgh(a_t) \approx a_t \bgh(a_t).
\]
Introduce the $\ns\times \na$ matrix $\Gh = [\bgh(a_t)]_{1\le t \le \na}$ with columns $\bgh(a_t)$ as well as the $\na\times \na$ diagonal matrix $\Lambda = \diag(a_t)$. The condition now reads
\[
M \Gh \approx \Gh \Lambda.
\]
When the columns of $\Gh$ are numerically linearly independent, we define the eigenmatrix for the real analytic case as
\[
M := \Gh \Lambda \Gh^+,
\]
where the pseudoinverse $\Gh^+$ is computed by thresholding the singular values of $\Gh$.

{\bf Combined with ESPRIT.}  Define the vector $\but=[\ut_j]_{1\le j\le \ns}$, where $\ut_j$ are the noisy observations. Consider the matrix $\bbm \but & M\but & \ldots & M^\L \but \ebm$ with $\L>\nx$, obtained from applying $M$ repeatitively. Since $\but \approx \sum_k \bg(x_k) w_k$ and $M\bg(x)\approx x \bg(x)$,
\[
\bbm  \but & M\but & \ldots & M^\L \but \ebm
\approx
\bbm \bg(x_1) & \ldots & \bg(x_\nx) \ebm
\bbm w_1 & &\\& \ddots & \\& & w_\nx\ebm
\bbm
1 &  x_1 & \ldots & (x_1)^\L \\
\vdots & \vdots & \ddots & \vdots \\
1 &  x_\nx & \ldots & (x_\nx)^\L
\ebm.
\]
Let $\tilde{U} \tilde{S}\tilde{V}^*$ be the rank-$\nx$ truncated SVD of this left-hand side. The matrix $\tilde{V}^*$ satisfies
\[
\tilde{V}^* \approx P
\bbm
1 &  x_1 & \ldots & (x_1)^\L \\
\vdots & \vdots & \ddots & \vdots \\
1 &  x_\nx & \ldots & (x_\nx)^\L
\ebm,
\]
where $P$ is an unknown non-degenerate $\nx \times \nx$ matrix. Let $\Zt_U$ and $\Zt_D$ be the submatrices obtained by excluding the first column and the last column, respectively, i.e.,
\[
\Zt_U \approx P
\bbm
1      & \ldots & (x_1)^{\L-1} \\
\vdots & \ddots & \vdots \\
1      & \ldots & (x_\nx)^{\L-1}
\ebm,
\quad
\Zt_D \approx P
\bbm
x_1 & \ldots & (x_1)^\L \\
\vdots & \ddots & \vdots \\
x_\nx & \ldots & (x_\nx)^\L
\ebm.
\]
By forming
\[
\Zt_D (\Zt_U)^+ \approx
P
\bbm
x_1 &  & \\
& \ddots & \\
& & x_\nx
\ebm
P^{-1},
\]
one obtains the estimates $\{\xt_k\}$ for $\{x_k\}$ by computing the eigenvalues of $\Zt_D (\Zt_U)^+$.

With $\{\xt_k\}$ available, the least square solve
\[
\min_{\wt_k} \sum_j \left|\sum_k G(s_j,\xt_k) \wt_k - \ut_j\right|^2 
\]
gives the estimators $\{\wt_k\}$ for $\{w_k\}$.

\section{Multidimensional case}\label{sec:method}

This section extends the eigenmatrix approach to the multidimensional case. We start with the 2D case and then move on to the higher dimensions. The main reason is that the 2D real case can be addressed with complex variables.

\subsection{2D case}
The key idea is to reduce the 2D real case to the 1D complex case (see for example \cite{fan2023efficient}). To simplify the discussion, assume first that $X$ is $[-1,1]^2$. Define for each $x=(x^1,x^2)$ the vector $\bg(x):=[G(s_j,x)]_{1\le j\le \ns}$ and $\gamma(x) := x^1 + i x^2\in\C$. The first step is to construct $M$ such that $M \bg(x) \approx \gamma(x) \bg(x)$ for $x\in [-1,1]^2$. Numerically, it is more robust to use the normalized vector $\bgh(x) = \bg(x)/\|\bg(x)\|$. The condition then becomes
\[
M \bgh(x) \approx \gamma(x) \bgh(x), \quad x\in [-1,1]^2.
\]
We enforce this condition on a two-dimensional Chebyshev grid $\{a_t\}_{1\le t \le \na}$ of size $T$ (i.e., with $\sqrt{T}$ points in each dimension) on the square, 
\[
M \bgh(a_t) \approx \gamma(a_t) \bgh(a_t).
\]
Define the $\ns\times \na$ matrix $\Gh = [\bgh(a_t)]_{1\le t \le \na}$ with $\bgh(a_t)$ as columns and also the $\na\times \na$ diagonal matrix $\Lambda = \diag(\gamma(a_t))$. The above condition can be written as
\[
M \Gh \approx \Gh \Lambda.
\]
When the columns of $\Gh$ are numerically linearly independent, this suggests 
\[
M := \Gh \Lambda \Gh^+,
\]
where the pseudoinverse $\Gh^+$ is computed by thresholding the singular values of $\Gh$.

\begin{remark}
  We claim that, for real analytic kernels $G(s,x)$, enforcing the condition at the Chebyshev grid $\{a_t\}$ is sufficient. To see this,
  \[
  M \bg(x) \approx M\left(\sum_t c_t(x) \bg(a_t)\right) = \sum_t c_t(x) M \bg(a_t) \approx \sum_t c_t(x) (\gamma(a_t) \bg(a_t)) \approx \gamma(x) \bg(x),
  \]
  where $c_t(x)$ is the Chebyshev quadrature for $x$ associated with grid $\{a_t\}$. Here, the first and third approximations use the convergence property of the Chebyshev quadrature for the analytic functions $\bg(x)$ and $\gamma(x)\bg(x)$, and the second approximation directly comes from $M \bg(a_t) \approx \gamma(a_t) \bg(a_t)$.
\end{remark}

Let $\but=[\ut_j]_{1\le j\le \ns}$ and form the following matrix by applying $M$ repeatitively to $\but$
\[
\bbm  \but & M\but & \ldots & M^\L \but \ebm
\approx
\bbm \bg(x_1) & \ldots & \bg(x_\nx) \ebm
\bbm w_1 & &\\& \ddots & \\& & w_\nx\ebm
\bbm
1 &  \gamma(x_1) & \ldots & \gamma(x_1)^\L \\
\vdots & \vdots & \ddots & \vdots \\
1 &  \gamma(x_\nx) & \ldots & \gamma(x_\nx)^\L
\ebm
\]
with $\L>\nx$. From the SVD, one can get estimates for $\{\gamma(x_k)\}$. Their real and imaginary parts give the approximate spikes $\{\xt_k\}$. Finally, the least square solve
\[
\min_{\wt_k} \sum_j \left|\sum_k G(s_j,\xt_k) \wt_k - \ut_j\right|^2 
\]
gives the estimators $\{\wt_k\}$ for $\{w_k\}$.

\begin{remark}
Having the columns of $\Gh$ to be numerically linearly independent is essential. For example, the method fails on the kernels such as $G(s,x) = \ln |s-x|$, which is Green's function of the Laplace equation and violates the linear independence due to the mean value theorem. The following argument explains why. Suppose for now that for any $x\in[-1,1]^2$, $M \bg(x) \approx \gamma(x) \bg(x)$.  Since $G(s,x)= \ln |s-x|$ satisfies the Laplace equation in $x$, the mean value theorem states that
\[
\bg(x) = \fint \bg(b) db,
\]
where the averaged integral is taken over a finite circle centered at $x$. Applying $M$ to both sides leads to
\begin{equation}
\gamma(x)\bg(x) \approx M\bg(x) = \fint M \bg(b) db \approx \fint \gamma(b) \bg(b) db.
\label{eq:false}
\end{equation}
However, the real and imaginary parts of $\gamma(x)\bg(x)$, i.e., $x^1 \bg(x)$ and $x^2 \bg(x)$, are not harmonic functions in $x$. Therefore, \eqref{eq:false} cannot hold for $G(s,x)= \ln |s-x|$ without introducing a finite error. Therefore, it is not possible to construct $M$ such that $M \bg(x) \approx \gamma(x) \bg(x)$ holds for all $x\in [-1,1]^2$. This argument also works for any kernel $G(s,x)$ that is the Green's functions of any linear elliptic PDE.
\end{remark}

\subsection{3D case}

To simplify the discussion, assume that $X$ is $[-1,1]^3$. Define for each $x$ the vector $\bg(x):=[G(s_j,x)]_{1\le j\le \ns}$. The first step is to construct three eigenmatrices $M^1, M^2, M^3$ such that for $x=(x^1,x^2,x^3) \in [-1,1]^3$.
\[
M^1 \bg(x) \approx x^1 \bg(x), \quad
M^2 \bg(x) \approx x^2 \bg(x), \quad
M^3 \bg(x) \approx x^3 \bg(x).
\]
Introducing a product Chebyshev grid $\{a_t=(a_t^1,a_t^2,a_t^3)\}_{1\le t \le \na}$ of size $\na$, i.e., with $\na^{1/3}$ points per dimension, we enforce the condition on this Chebyshev grid
\[
M^1 \bgh(a_t) \approx a_t^1 \bgh(a_t), \quad
M^2 \bgh(a_t) \approx a_t^2 \bgh(a_t), \quad
M^3 \bgh(a_t) \approx a_t^3 \bgh(a_t).
\]
Define the $\ns\times \na$ matrix $\Gh = [\bgh(a_t)]_{1\le t \le \na}$ with $\bgh(a_t)$ as columns and also $\na\times \na$ diagonal matrices $\Lambda^1 = \diag(a_t^1)$, $\Lambda^2 = \diag(a_t^2)$, and $\Lambda^3 = \diag(a_t^3)$. The above condition can be written in a matrix form as
\[
M^1 \Gh \approx \Gh \Lambda^1, \quad
M^2 \Gh \approx \Gh \Lambda^2, \quad
M^3 \Gh \approx \Gh \Lambda^3.
\]
When the columns of $\Gh$ are numerically linearly independent, this suggests the following choice of the eigenmatrix,
\[
M^1 := \Gh \Lambda^1 \Gh^+, \quad
M^2 := \Gh \Lambda^2 \Gh^+, \quad
M^3 := \Gh \Lambda^3 \Gh^+.
\]
where the pseudoinverse $\Gh^+$ is computed by thresholding the singular values of $\Gh$.

\begin{remark}
  The matrices $M^1$, $M^2$, and $M^3$ approximately commute since they share the same set of approximate eigenvalues and eigenvectors. Therefore, any product formed from them is approximately independent of the order.
\end{remark}

\begin{remark}
  Enforcing the condition at the Chebyshev grid $\{a_t\}$ is again sufficient, i.e., $M^1 \bg(a_t) \approx a_t^1 \bg(a_t)$ implies $M^1 \bg(x) \approx x^1 \bg(x)$ for all $x\in [-1,1]^3$. 
  \[
  M^1 \bg(x) \approx M^1\left(\sum_t c_t(x) \bg(a_t)\right) = \sum_t c_t(x) M^1 \bg(a_t) \approx \sum_t c_t(x) (a_t^1 \bg(a_t)) \approx x^1 \bg(x),
  \]
  where $c_t(x)$ is the Chebyshev quadrature for $x$ associated with grid $\{a_t\}$. Here, the first and third approximations use the convergence property of the Chebyshev quadrature for analytic functions $\bg(x)$ and $x^1\bg(x)$, and the second approximation directly comes from $M^1 \bg(a_t) \approx a_t^1 \bg(a_t)$. The same derivation holds for $M^2$ and $M^3$.
\end{remark}

With the eigenmatrices available, the rest of the algorithm follows the work in \cite{andersson2018esprit}. For any triple $\alpha=(\alpha^1,\alpha^2,\alpha^3)$, we define $(x)^\alpha := (x^1)^{\alpha^1} (x^2)^{\alpha^2} (x^3)^{\alpha^3}$ and $M^\alpha := (M^1)^{\alpha^1} (M^2)^{\alpha^2} (M^3)^{\alpha^3}$. Since $M^1$, $M^2$, and $M^3$ approximately commute, the definition of $M^\alpha$ is approximately independent of the order. Consider the matrix
\[
\bbm  \but &\ldots& M^\alpha\but &\ldots& M^{(\L,\L,\L)}\but \ebm
\]
with $(\L+1)^3$ columns, ordered in the lexicographical order. We have the following approximation
\begin{align*}
        & \bbm \but & \ldots & M^\alpha\but & \ldots &  M^{(\L,\L,\L)}\but\ebm \\
\approx &
\bbm \bg(x_1) & \ldots & \bg(x_\nx) \ebm
\bbm w_1 & &\\& \ddots & \\& & w_\nx\ebm
\bbm
1      & \ldots & (x_1)^\alpha & \ldots & (x_1)^{(\L,\L,\L)}\\
\vdots & \ddots & \vdots & \ddots & \vdots \\
1      & \ldots & (x_\nx)^\alpha & \ldots & (x_\nx)^{(\L,\L,\L)}
\ebm.
\end{align*}
Notice that one never needs to form the matrices $M^\alpha$ explicitly. The columns $M^\alpha \but$ can be computed effectively by progressively applying only the matrices $M^1$, $M^2$, and $M^3$.

Let $\tilde{U} \tilde{S}\tilde{V}^*$ be the rank-$\nx$ truncated SVD of this matrix. The matrix $\tilde{V}^*$ satisfies
\[
\tilde{V}^* \approx P
\bbm
1 & \ldots & (x_1)^\alpha & \ldots & (x_1)^{(\L,\L,\L)}\\
\vdots & \ddots & \vdots & \ddots & \vdots\\
1 & \ldots & (x_\nx)^\alpha & \ldots & (x_\nx)^{(\L,\L,\L)}
\ebm,
\]
where $P$ is an unknown non-degenerate $\nx \times \nx$ matrix.

For the first dimension, let $\Zt^1_U$ be the submatrix that excludes the columns $\alpha$ with $\alpha^1=\L$ and
$\Zt^1_D$ be the submatrix that excludes the columns $\alpha$ with $\alpha^1=0$. Then 
\[
\Zt^1_U \approx P
\bbm
1      & \ldots & (x_1)^{(\L-1,\L,\L)} \\
\vdots & \ddots & \vdots \\
1      & \ldots & (x_\nx)^{(\L-1,\L,\L)}
\ebm,
\quad
\Zt^1_D \approx P
\bbm
(x_1)^{(1,0,0)} & \ldots & (x_1)^{(\L,\L,\L)} \\
\vdots & \ddots & \vdots \\
(x_\nx)^{(1,0,0)} & \ldots & (x_\nx)^{(\L,\L,\L)}
\ebm.
\]
Introduce
\[
B^1 := \Zt^1_D (\Zt^1_U)^+ \approx
P
\bbm
x_1^1 &  & \\
& \ddots & \\
& & x_\nx^1
\ebm
P^{-1}.
\]

For the second dimension, let $\Zt^2_U$ be the submatrix that excludes the columns $\alpha$ with $\alpha^2=\L$ and
$\Zt^2_D$ be the submatrix that excludes the columns $\alpha$ with $\alpha^2=0$. Then 
\[
\Zt^2_U \approx P
\bbm
1      & \ldots & (x_1)^{(\L,\L-1,\L)} \\
\vdots & \ddots & \vdots \\
1      & \ldots & (x_\nx)^{(\L,\L-1,\L)}
\ebm,
\quad
\Zt^2_D \approx P
\bbm
(x_1)^{(0,1,0)} & \ldots & (x_1)^{(\L,\L,\L)} \\
\vdots & \ddots & \vdots \\
(x_\nx)^{(0,1,0)} & \ldots & (x_\nx)^{(\L,\L,\L)}
\ebm.
\]
Introduce
\[
B^2 := \Zt^2_D (\Zt^2_U)^+ \approx
P
\bbm
x_1^2 &  & \\
& \ddots & \\
& & x_\nx^2
\ebm
P^{-1}.
\]

For the third dimension, let $\Zt^3_U$ be the submatrix that excludes the columns $\alpha$ with $\alpha^3=\L$ and
$\Zt^3_D$ be the submatrix that excludes the columns $\alpha$ with $\alpha^3=0$. Then 
\[
\Zt^3_U \approx P
\bbm
1      & \ldots & (x_1)^{(\L,\L,\L-1)} \\
\vdots & \ddots & \vdots \\
1      & \ldots & (x_\nx)^{(\L,\L,\L-1)}
\ebm,
\quad
\Zt^2_D \approx P
\bbm
(x_1)^{(0,0,1)} & \ldots & (x_1)^{(\L,\L,\L)} \\
\vdots & \ddots & \vdots \\
(x_\nx)^{(0,0,1)} & \ldots & (x_\nx)^{(\L,\L,\L)}
\ebm.
\]
Introduce
\[
B^3 := \Zt^3_D (\Zt^3_U)^+ \approx
P
\bbm
x_1^3 &  & \\
& \ddots & \\
& & x_\nx^3
\ebm
P^{-1}.
\]

The definitions of $B^1$, $B^2$, and $B^3$ suggest that their eigenvalues approximate the three components of $\{x_k\}$, respectively. In order to avoid matching these components, it is more convenient to recover $P$ first. However, when $B^1$, $B^2$, or $B^3$ have degenerate eigenvalues, the individual eigenvalue decompositions might give incorrect answers. To robustly recover $P$, we choose three random unit complex numbers $\beta^1,\beta^2,\beta^3$ and define
\[
B := \beta^1 B^1 + \beta^2 B^2 + \beta^3 B^3 \approx
P
\bbm
\beta^1 x_1^1 + \beta^2 x_1^2 + \beta^3 x_1^3 &  & \\
& \ddots & \\
& & \beta^1 x_\nx^1 + \beta^2 x_\nx^2 + \beta^3 x_\nx^3
\ebm
P^{-1}.
\]
With high probability, these eigenvalues are well-separated from each other. Therefore, computing the eigenvectors of $B$ gives $P$. From $P$ and $P^{-1}$, we can approximate $x_1^1,\ldots,x_\nx^1$ from $B^1$, $x_1^2,\ldots,x_\nx^2$ from $B^2$, and $x_1^3,\ldots,x_\nx^3$ from $B^3$, respectively. Combining them gives the approximation to $x_1,\ldots, x_\nx$. Finally, the least square solve
\[
\min_{\wt_k} \sum_j \left|\sum_k G(s_j,\xt_k) \wt_k - \ut_j\right|^2 
\]
leads to the estimators $\{\wt_k\}$ for $\{w_k\}$.

\subsection{Higher dimensional real and complex cases}
For the $d$-dimensional real case, assume that the domain $X$ is $[-1,1]^d$. We can choose a product Chebyshev grid to construct the eigenmatrices $M^1,\ldots, M^d$, each responsible for
\[
M^t \bg(x) \approx x^t \bg(x), \quad 1\le t \le d.
\]
A general multi-index power for $\alpha=(\alpha^1,\ldots,\alpha^d)$ is defined in the same way
\[
M^\alpha := (M^1)^{\alpha^1} \ldots (M^d)^{\alpha^d}.
\]
From the matrix
\[
\bbm  \but &\ldots& M^\alpha\but &\ldots& M^{(\L,\cdots,\L)}\but \ebm
\]
with $(L+1)^d$ columns, we compute $\Zt^t_U$ and $\Zt^t_D$ for each dimension $t=1,\ldots,d$. From the matrices $B^t := \Zt^t_D (\Zt^t_U)^+$ and $B := \beta^1 B^1 + \ldots + \beta^d B^d$, we can approximate the spikes $\{x_k\}$.

The higher dimensional complex case is also similar, where the only difference is to replace the product Chebyshev grid with a product uniform grid over the circles.

\section{Numerical results}\label{sec:num}

This section applies the eigenmatrix approach to the multidimensional unstructured sparse recovery problems mentioned above.

\subsection{2D case}

In the examples, the Chebyshev grid is $32\times 32$, i.e., $\na=32^2$. The spike weights $\{w_k\}$ are set to be $1$ and the noises $\{Z_j\}$ are Gaussian. In each plot, blue, green, and red stand for the exact solution, the result before postprocessing, and the one after postprocessing, respectively.

\begin{example}[Sparse deconvolution] The problem setup is 
  \begin{itemize}
  \item $G(s,x) = \frac{1}{\|s-x\|}$.
  \item $X=[-1,1]^2$.
  \item $\{s_j\}$ are $J=1024$ random points in $[-2,2]^2$ outside $X$.
  \end{itemize}
  
  \begin{figure}[h!]
    \centering
    \includegraphics[scale=0.25]{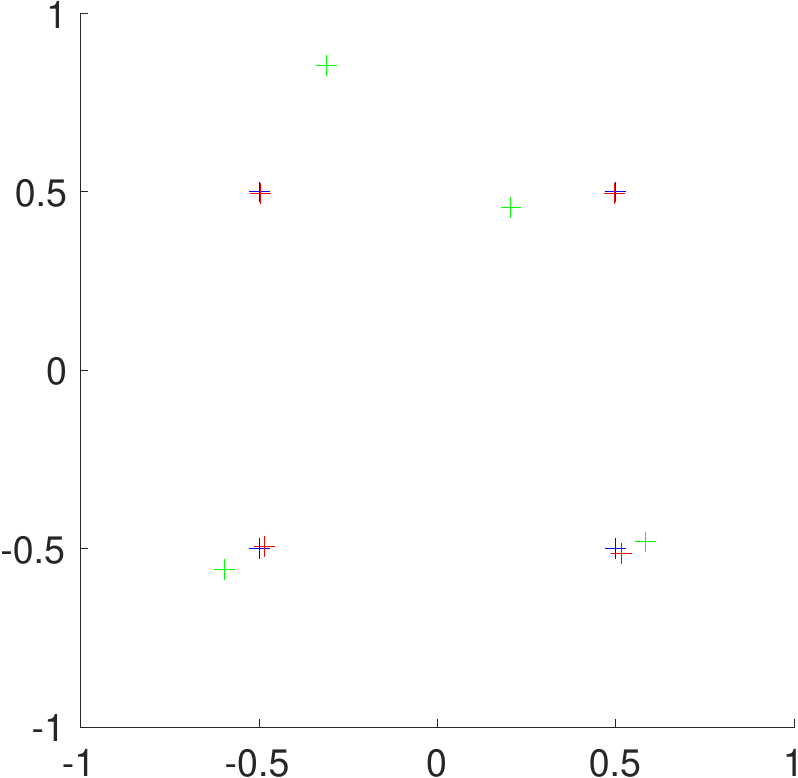}
    \includegraphics[scale=0.25]{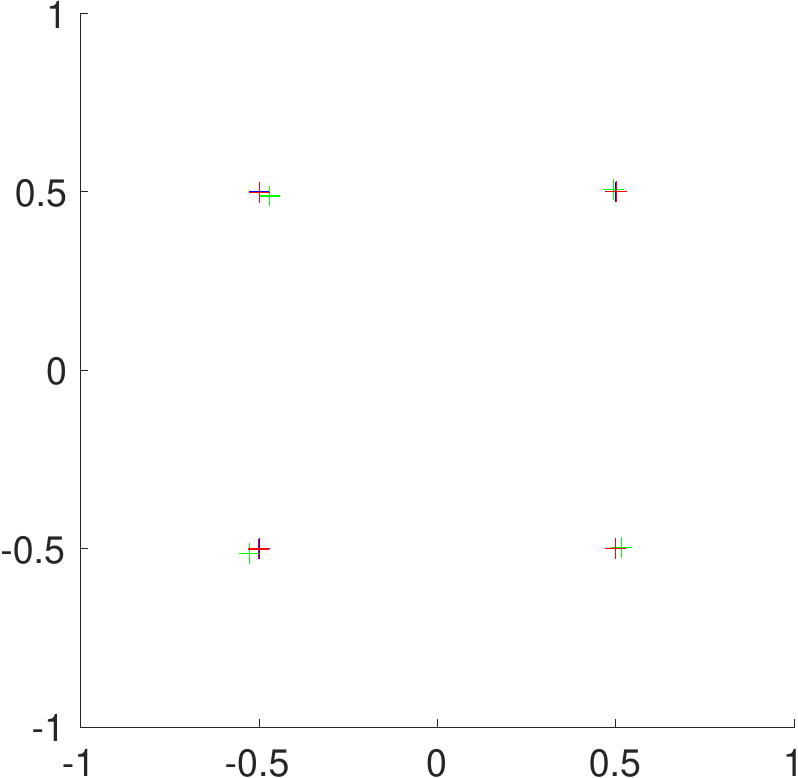}
    \includegraphics[scale=0.25]{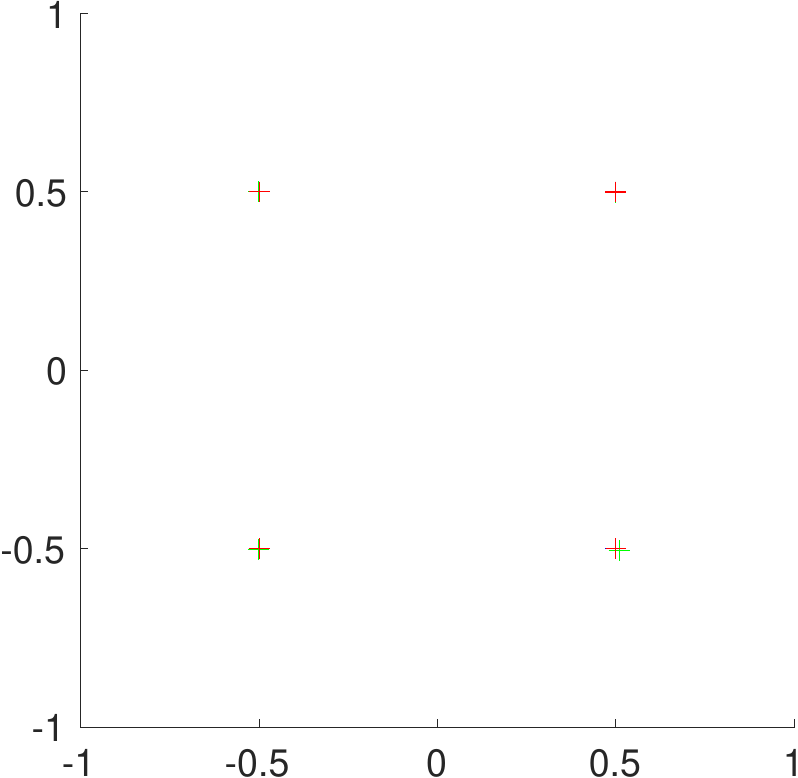}\\
    \includegraphics[scale=0.25]{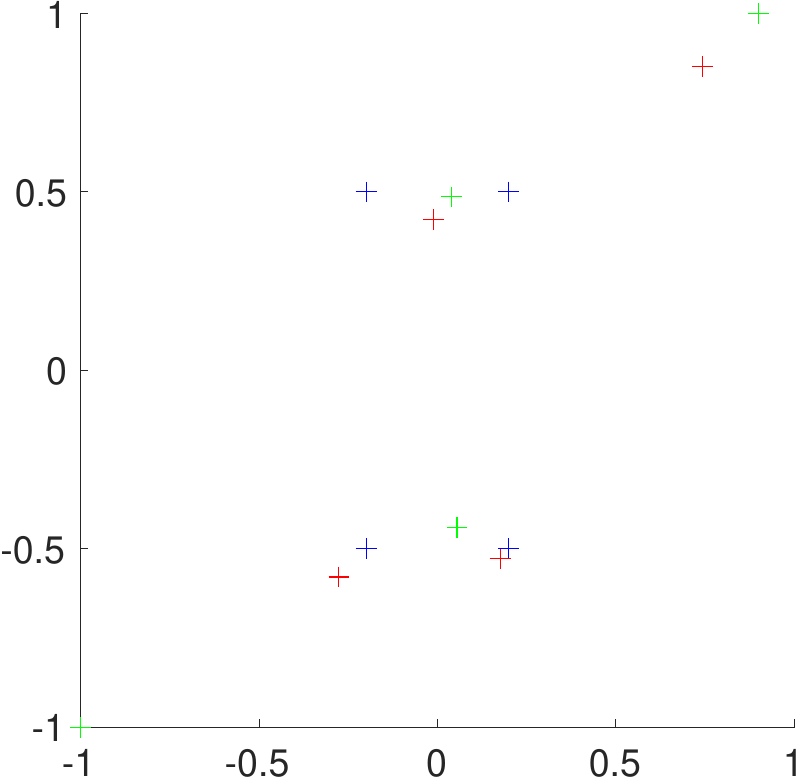}
    \includegraphics[scale=0.25]{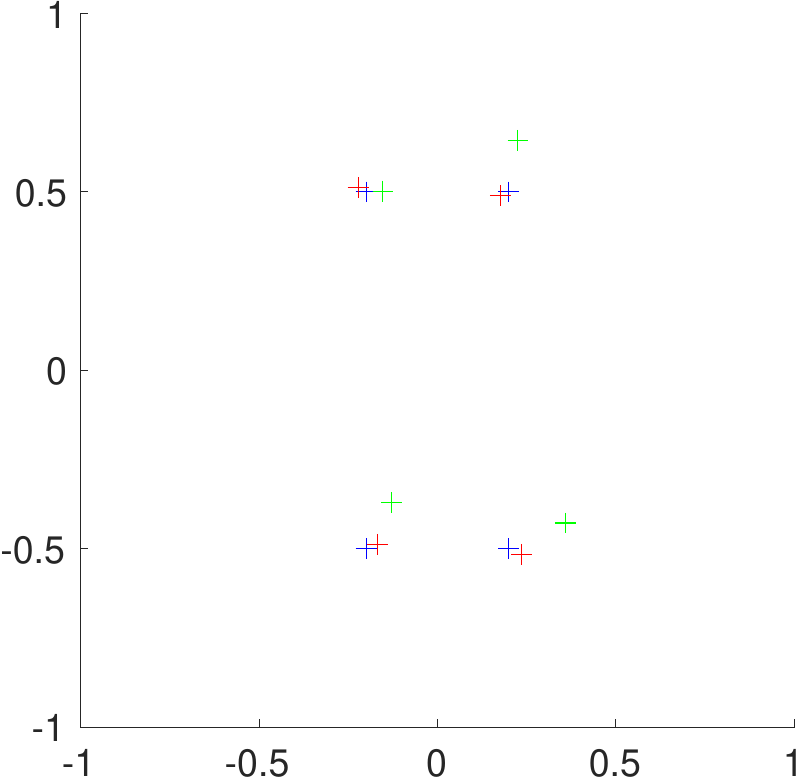}
    \includegraphics[scale=0.25]{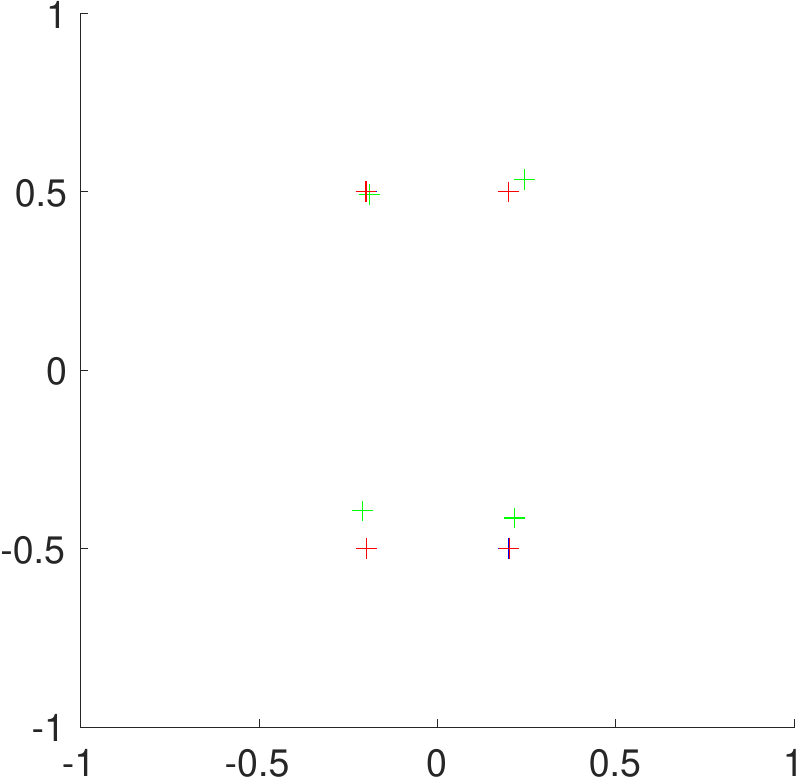} 
    \caption{Sparse deconvolution. $G(s,x) = \frac{1}{\|s-x\|}$.  $X=[-1,1]^2$. $\{s_j\}$ are $J=1024$ random points in $[-2,2]^2$ outside $X$. Columns: $\sigma$ equals to $10^{-2}$, $10^{-3}$, and $10^{-4}$, respectively. Rows: the easy case (top) and the hard case (bottom).}
    \label{fig:ex2P}
  \end{figure}
  
  Figure \ref{fig:ex2P} summarizes the experimental results. The three columns correspond to noise levels $\sigma$ equal to $10^{-2}$, $10^{-3}$, and $10^{-4}$. Two tests are performed. In the first one (top row), $\{x_k\}$ are well-separated from each other.  The plots show accurate recovery of the spikes from all $\sigma$ values. In the second one (bottom row), the spikes form two nearby pairs. The reconstruction at $\sigma=10^{-2}$ shows a noticeable error, while the results for $\sigma=10^{-3}$ and $\sigma=10^{-4}$ are accurate.
\end{example}

\begin{example}[Fourier inversion] The problem setup is
  \begin{itemize}
  \item $G(s,x) = \exp(\pi i s \cdot x)$.
  \item $X=[-1,1]^2$.
  \item $\{s_j\}$ are randomly chosen $J=1024$ points in $[-8,8]^2$.
  \end{itemize}

  \begin{figure}[h!]
  \centering
  \includegraphics[scale=0.25]{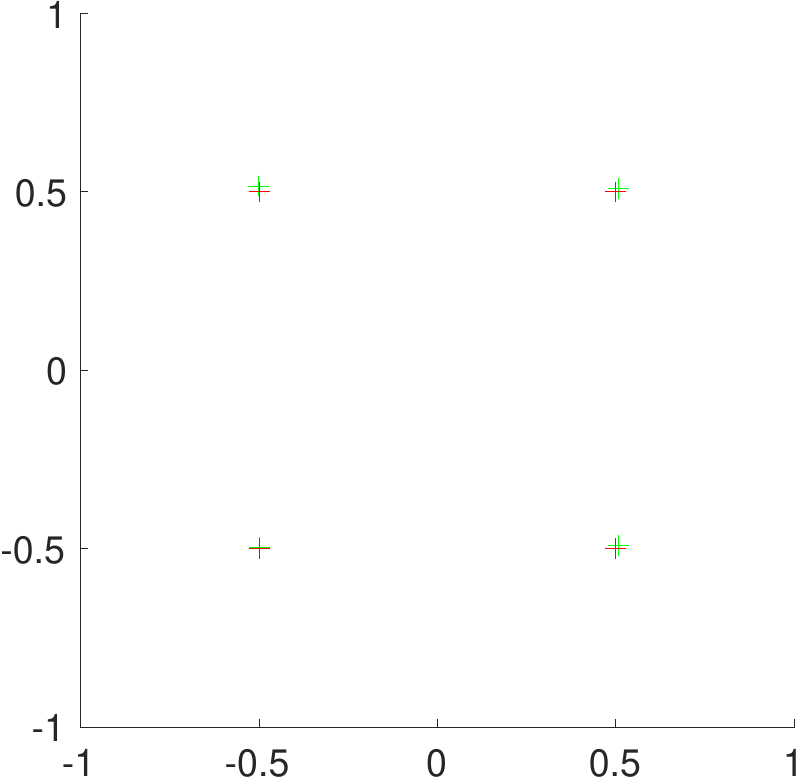}
  \includegraphics[scale=0.25]{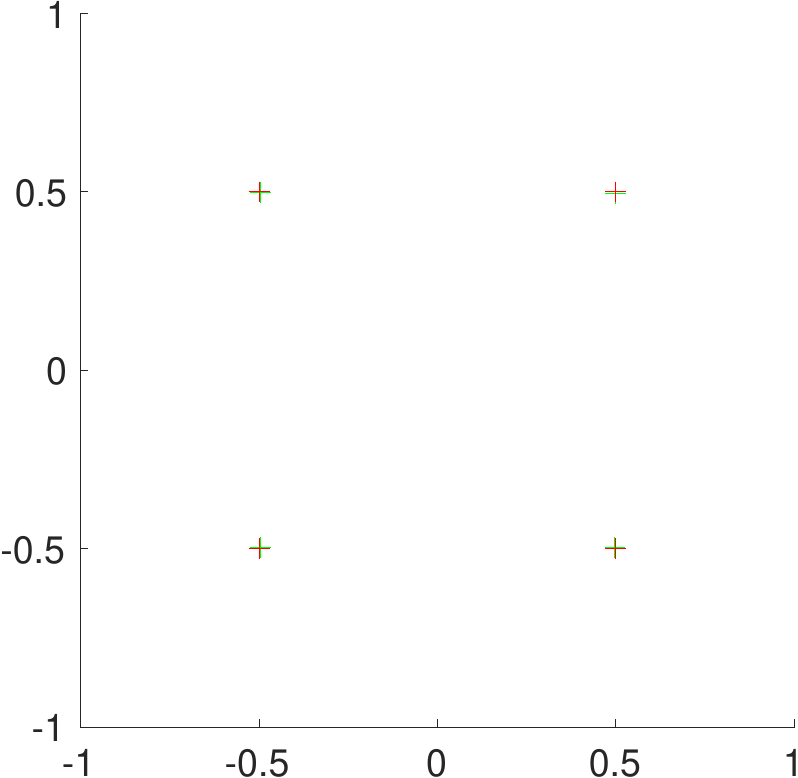}
  \includegraphics[scale=0.25]{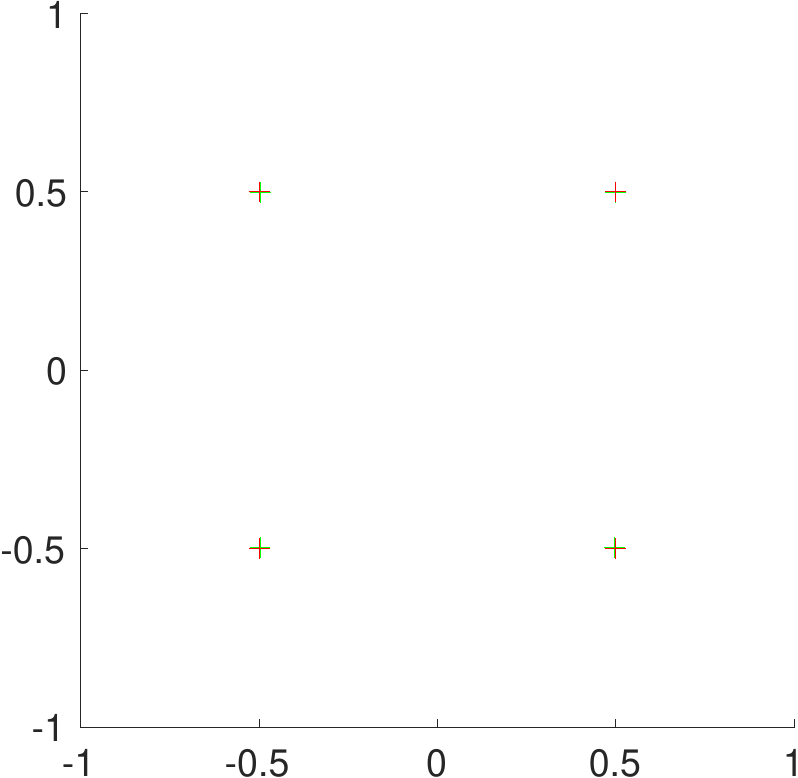}\\
  \includegraphics[scale=0.25]{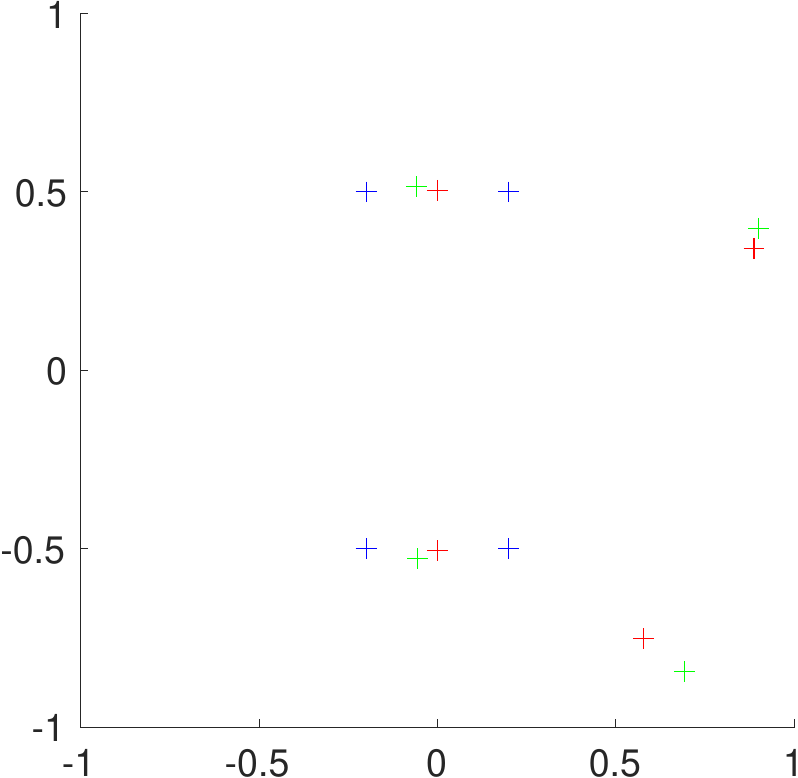}
  \includegraphics[scale=0.25]{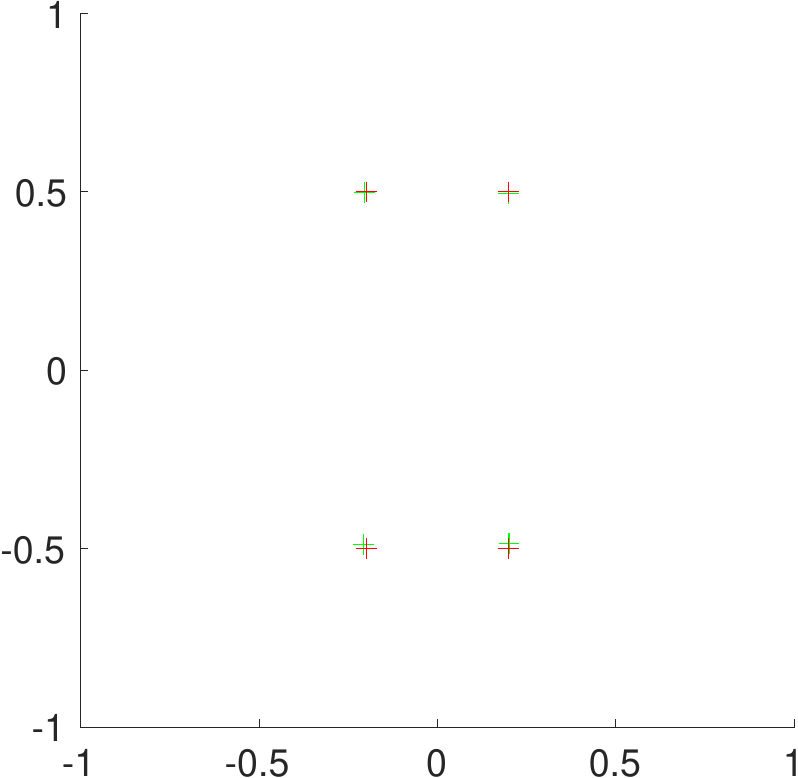}
  \includegraphics[scale=0.25]{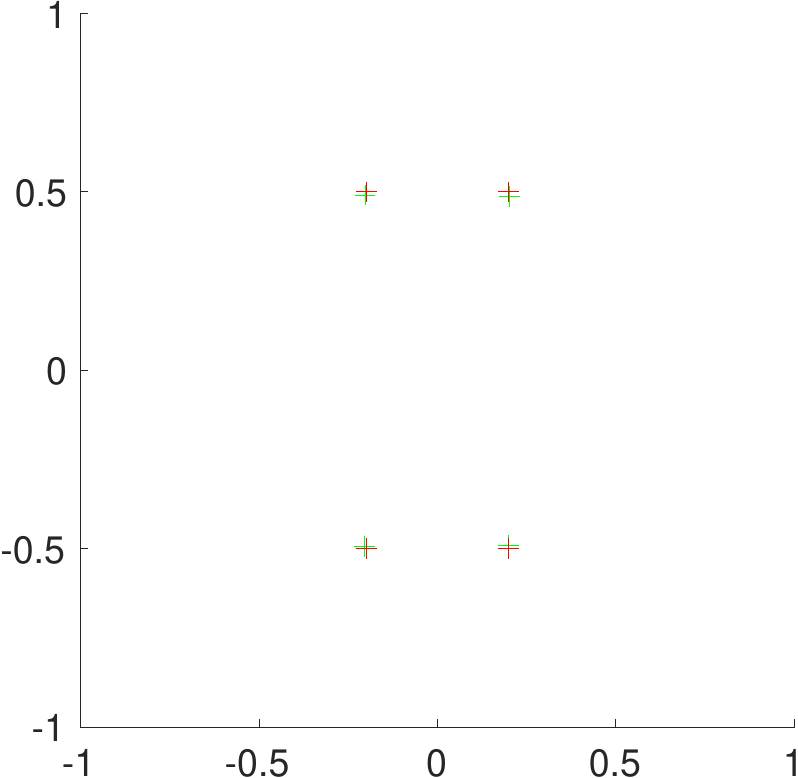} 
  \caption{Fourier inversion. $G(s,x) = \exp(\pi i s \cdot x)$. $X=[-1,1]^2$. $\{s_j\}$ are $J=1024$ randomly chosen points in $[-8,8]^2$. Columns: $\sigma$ equals to $10^{-2}$, $10^{-3}$, and $10^{-4}$, respectively. Rows: the easy case (top) and the hard case (bottom).}
  \label{fig:ex2F}
  \end{figure}

  Figure \ref{fig:ex2F} summarizes the experimental results. The three columns correspond to $\sigma$ equal to $10^{-2}$, $10^{-3}$, and $10^{-4}$. Two tests are performed. In the first one (top row), $\{x_k\}$ are well-separated, and the reconstructions are accurate for all $\sigma$ values. In the second one (bottom row), the spikes form two nearby pairs.  The reconstruction at $\sigma=10^{-2}$ shows a noticeable error, but the ones for lower noises are accurate.
\end{example}

\subsection{3D case}
In the 3D examples, the Chebyshev grid is $16\times 16 \times 16$, i.e., $\na=16^3$. The spike weights $\{w_k\}$ are set to be $1$ and the noises $\{Z_j\}$ are still Gaussian.

\begin{example}[Sparse deconvolution] The problem setup is 
  \begin{itemize}
  \item $G(s,x) = \frac{1}{\|s-x\|^{1/2}}$.
  \item $X=[-1,1]^3$.
  \item $\{s_j\}$ are $J=8192$ random points in $[-2,2]^3$ outside $X$.
  \end{itemize}

  \begin{figure}[h!]
    \centering
    \includegraphics[scale=0.25]{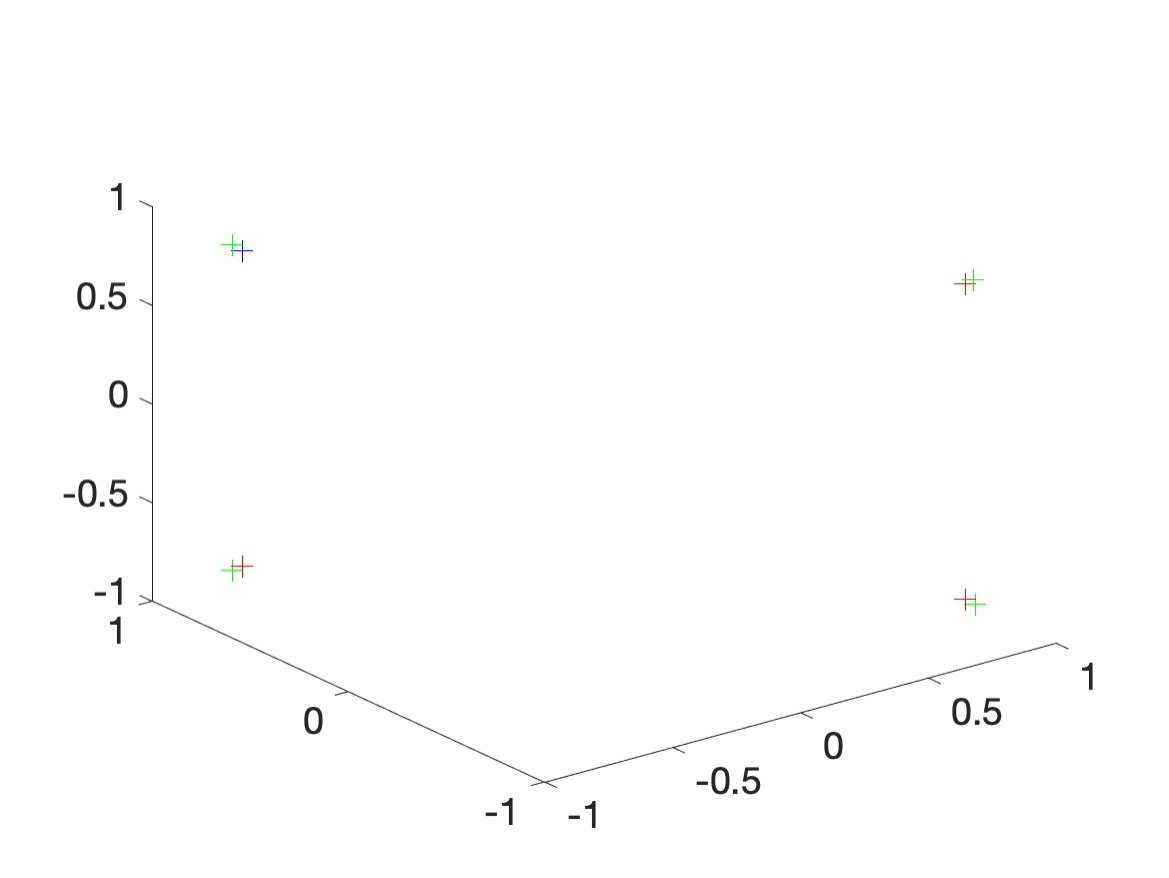}
    \includegraphics[scale=0.25]{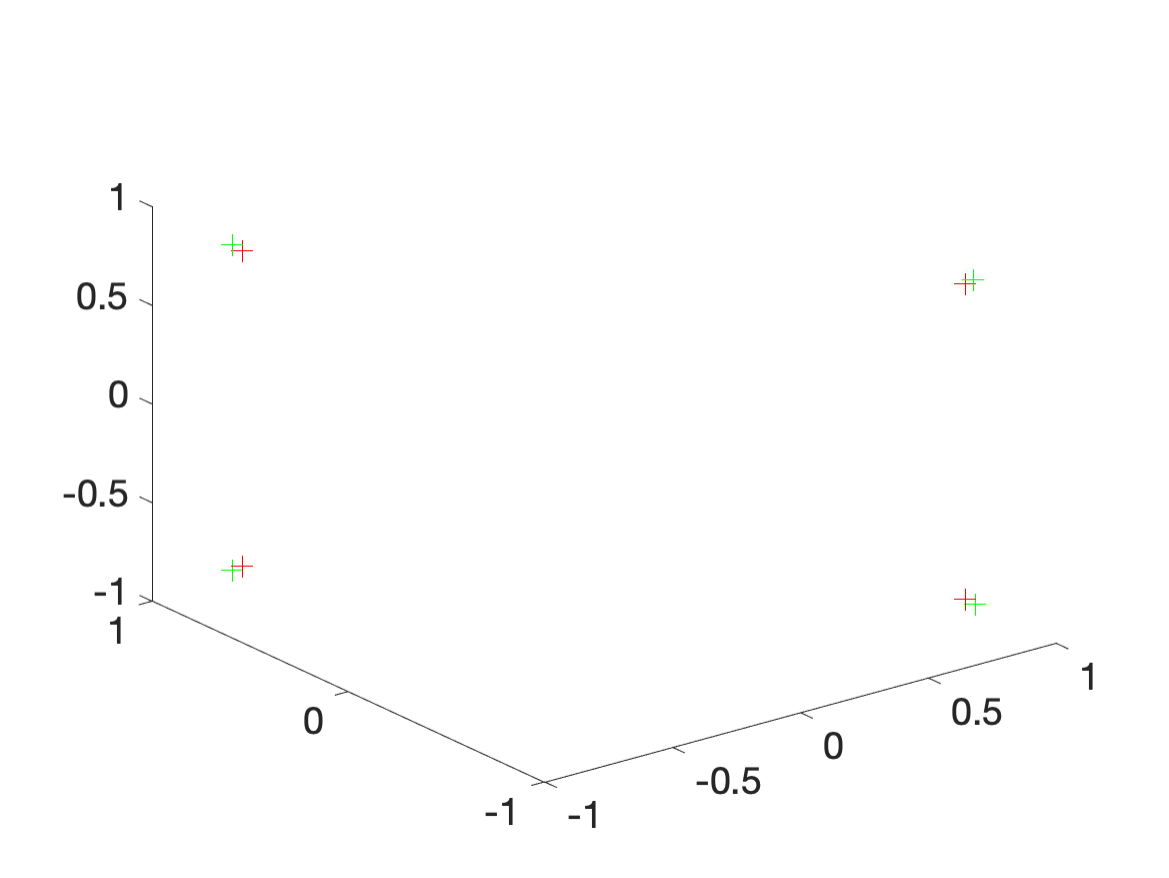}
    \includegraphics[scale=0.25]{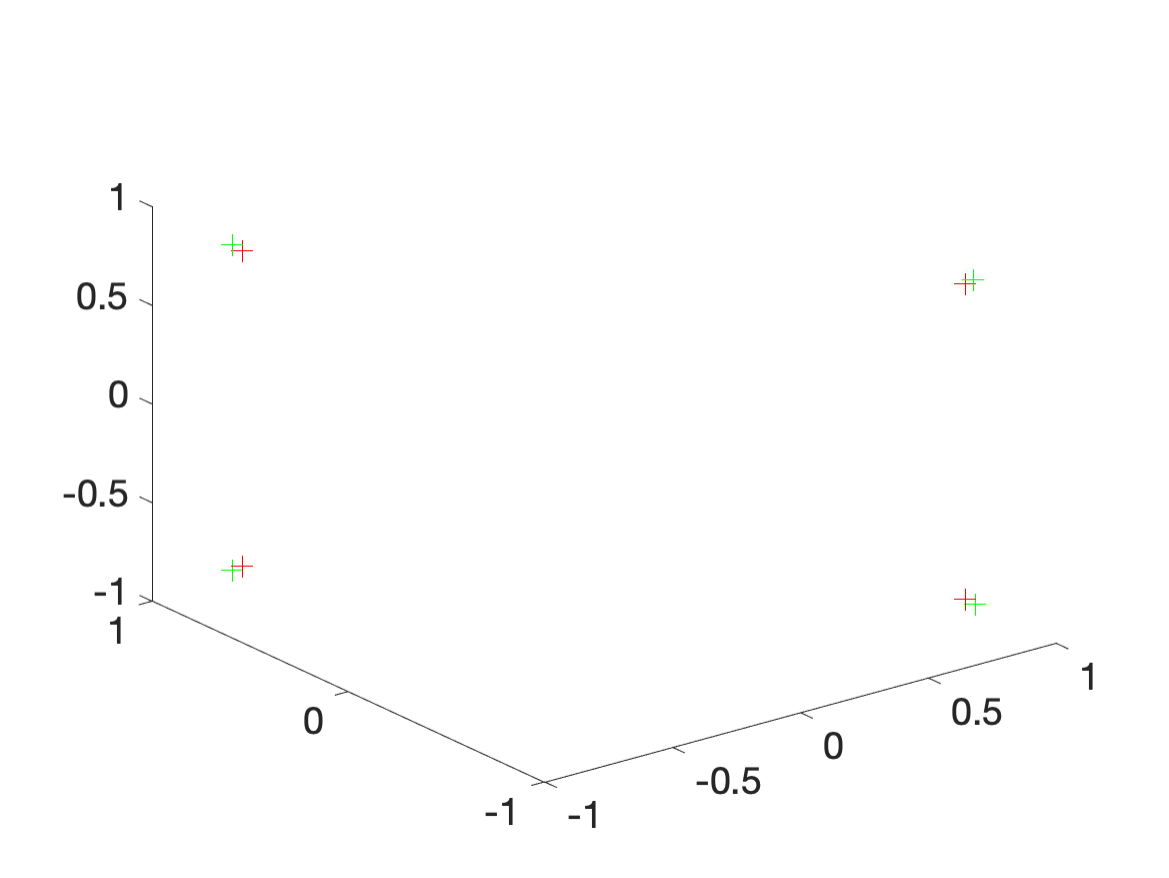}\\
    \includegraphics[scale=0.25]{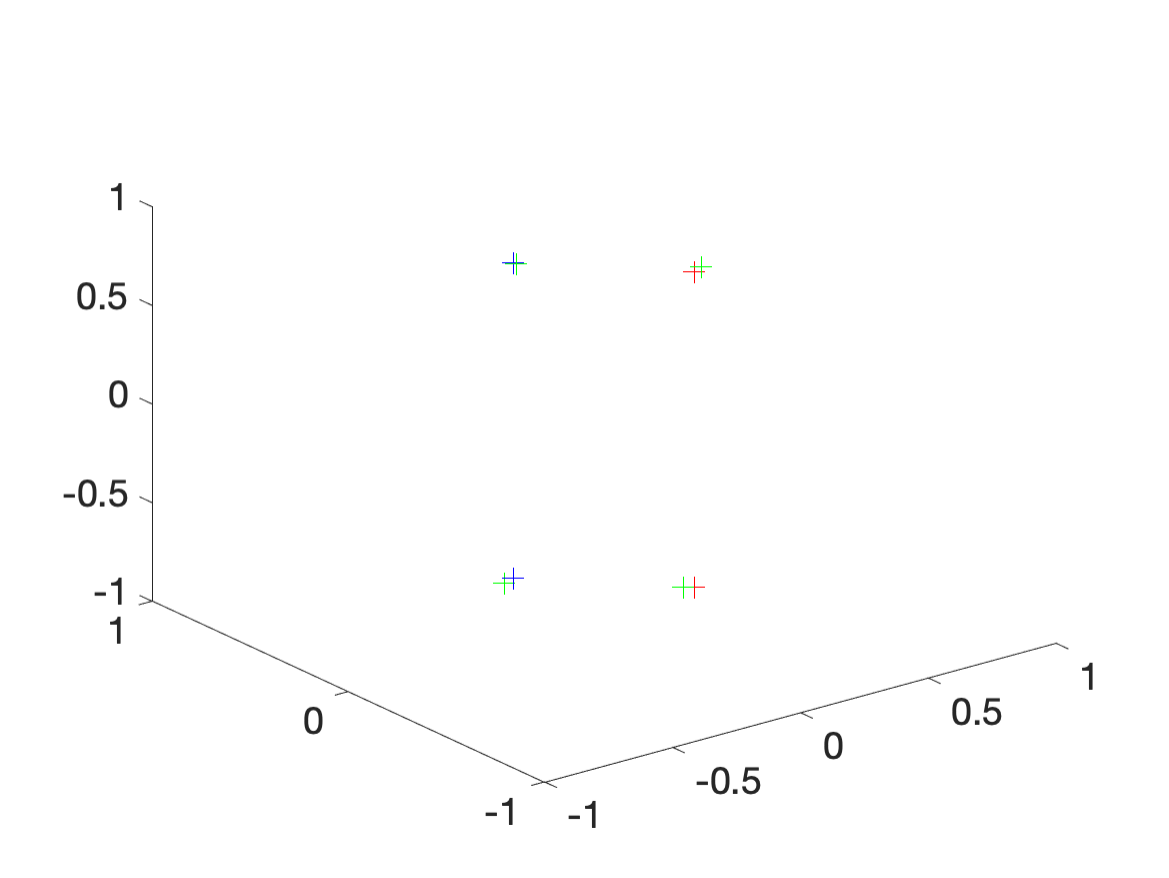}
    \includegraphics[scale=0.25]{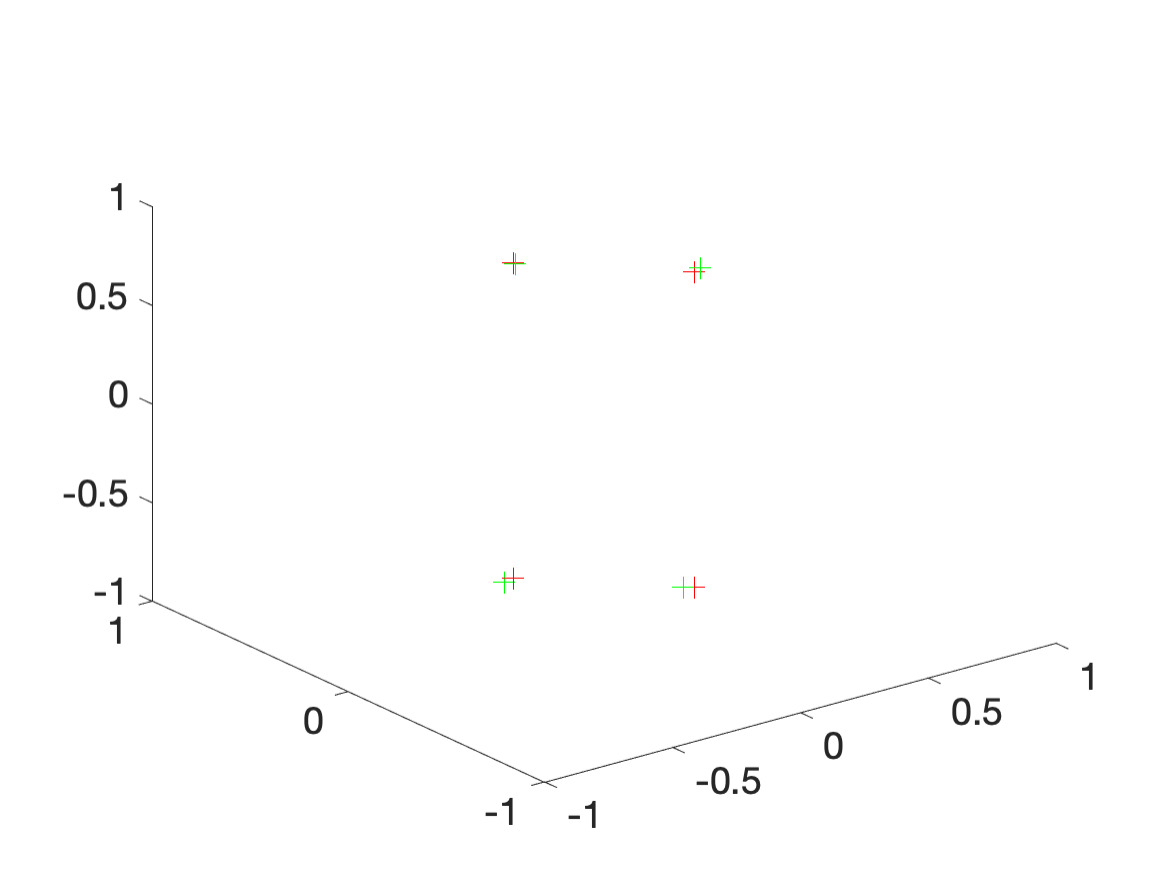}
    \includegraphics[scale=0.25]{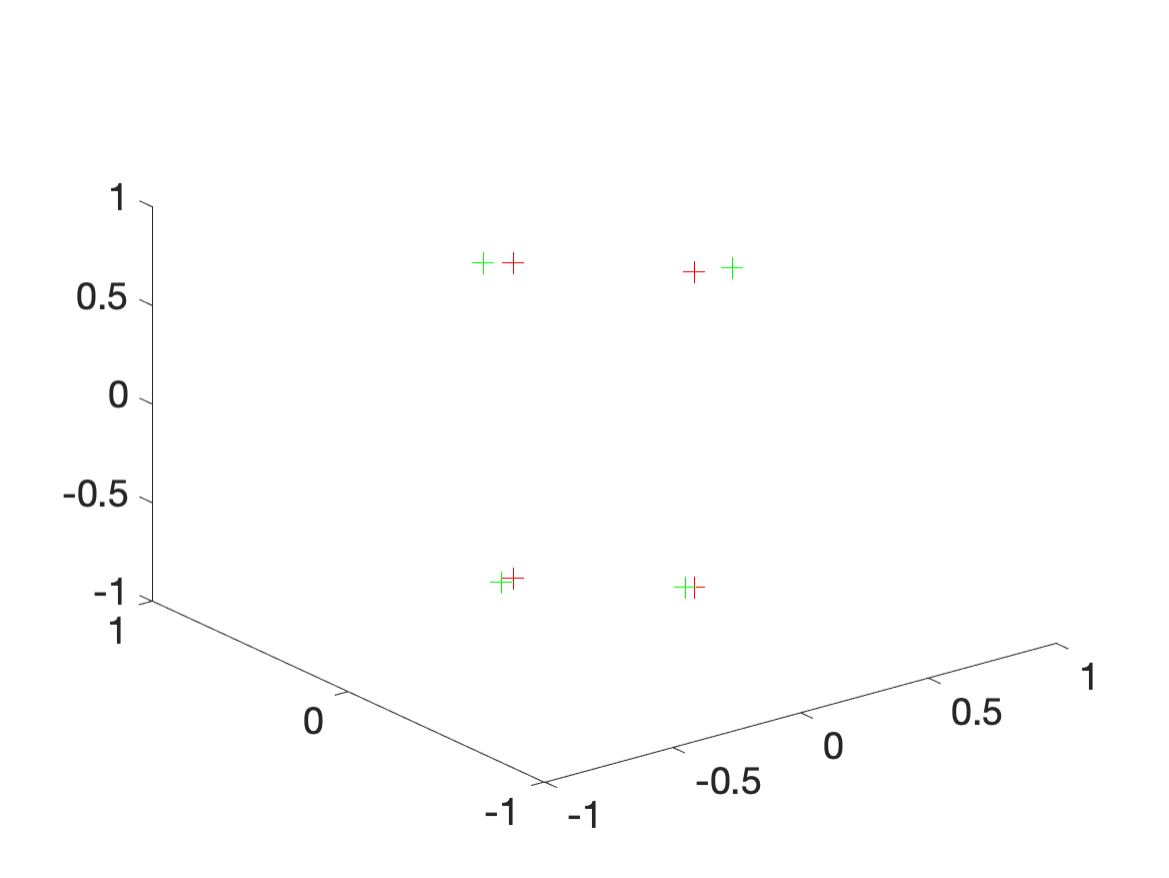} 
    \caption{Sparse deconvolution. $G(s,x) = \frac{1}{\|s-x\|^{1/2}}$.  $X=[-1,1]^3$. $\{s_j\}$ are $J=8192$ random points in $[-2,2]^3$ outside $X$. Columns: $\sigma$ equals to $10^{-4}$, $10^{-5}$, and $10^{-6}$, respectively. Rows: the easy case (top) and the hard case (bottom).}
    \label{fig:ex3P}
  \end{figure}
  
  Figure \ref{fig:ex3P} summarizes the experimental results. The three columns correspond to noise levels $\sigma$ equal to $10^{-4}$, $10^{-5}$, and $10^{-6}$. Two tests are performed. In the first one (top row), $\{x_k\}$ are well-separated from each other.  In the second one (bottom row), the spikes form two nearby pairs. The plots show accurate recovery of the spikes from all $\sigma$ values in both cases.
\end{example}

\begin{example}[Fourier inversion] The problem setup is
  \begin{itemize}
  \item $G(s,x) = \exp(\pi i s \cdot x)$.
  \item $X=[-1,1]^3$.
  \item $\{s_j\}$ are $J=8192$ randomly chosen points in $[-4,4]^3$.
  \end{itemize}

  \begin{figure}[h!]
  \centering
  \includegraphics[scale=0.25]{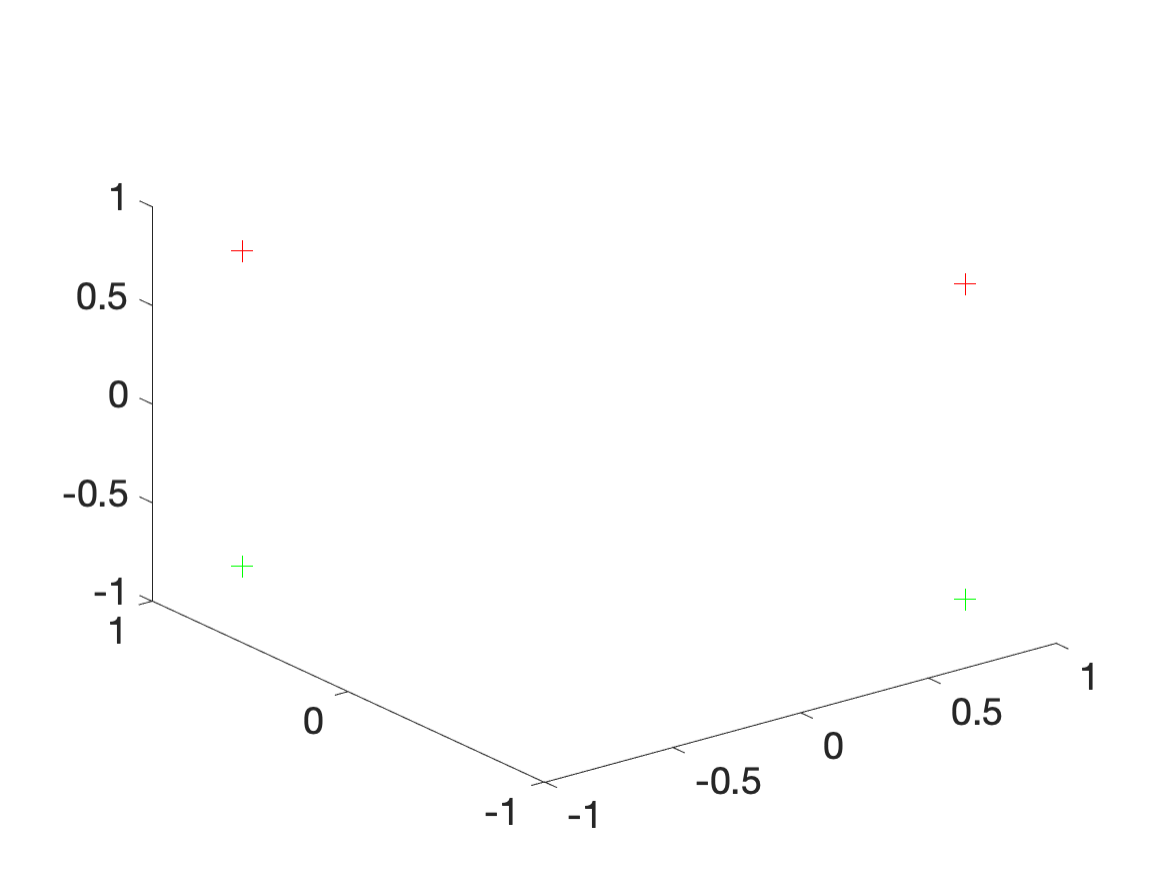}
  \includegraphics[scale=0.25]{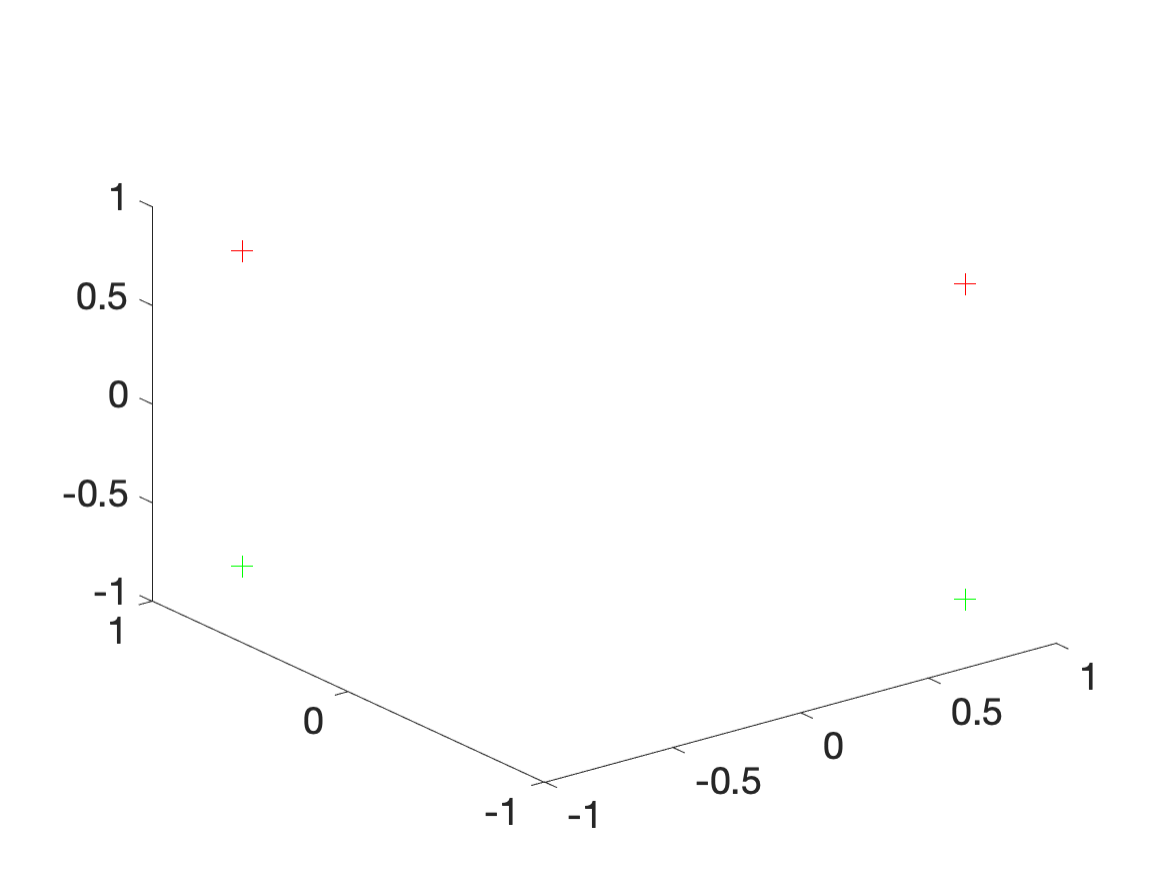}
  \includegraphics[scale=0.25]{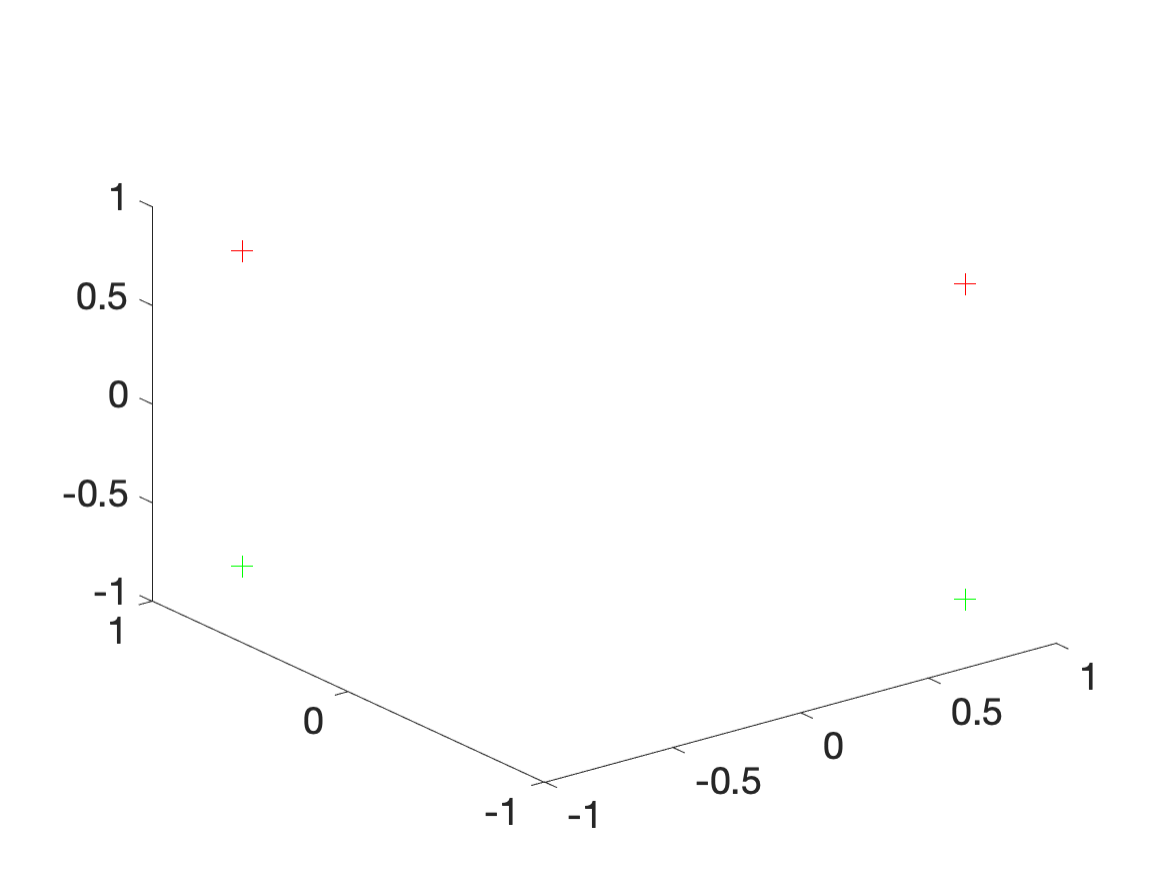}\\
  \includegraphics[scale=0.25]{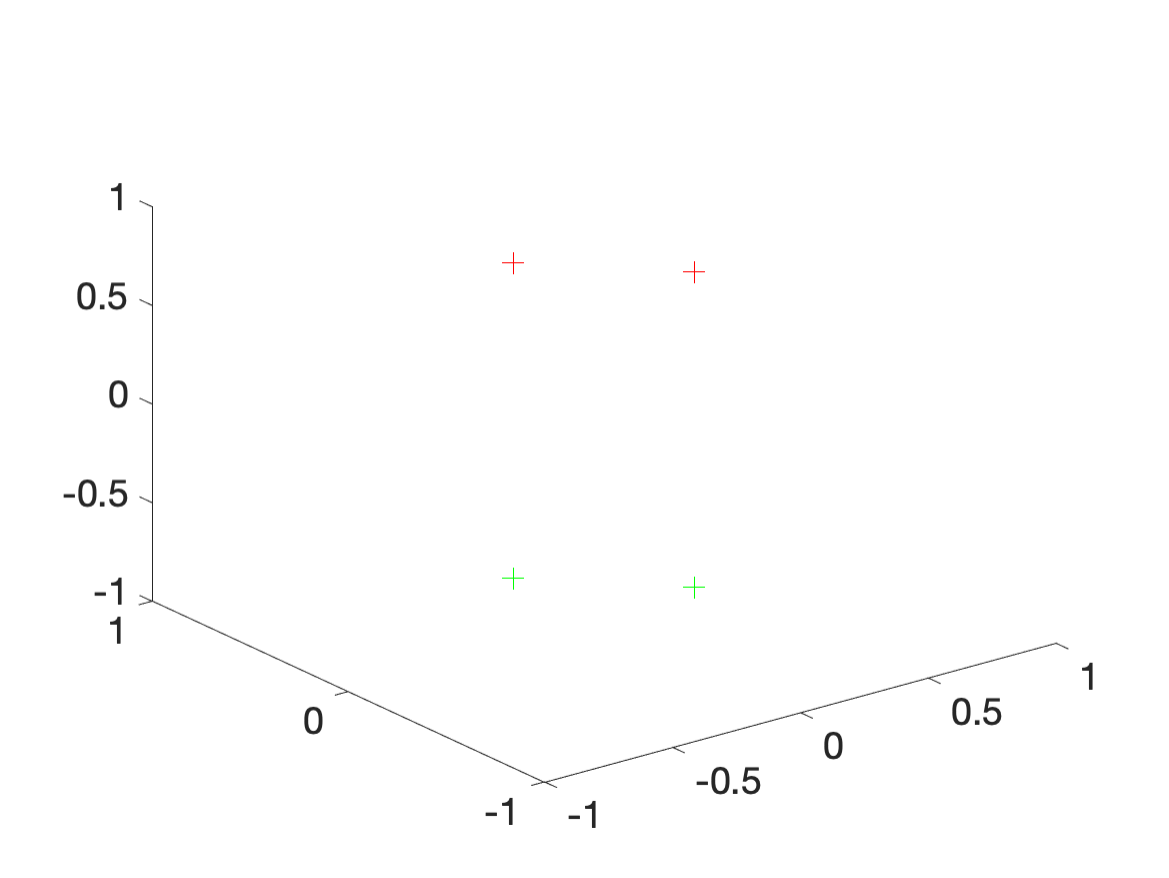}
  \includegraphics[scale=0.25]{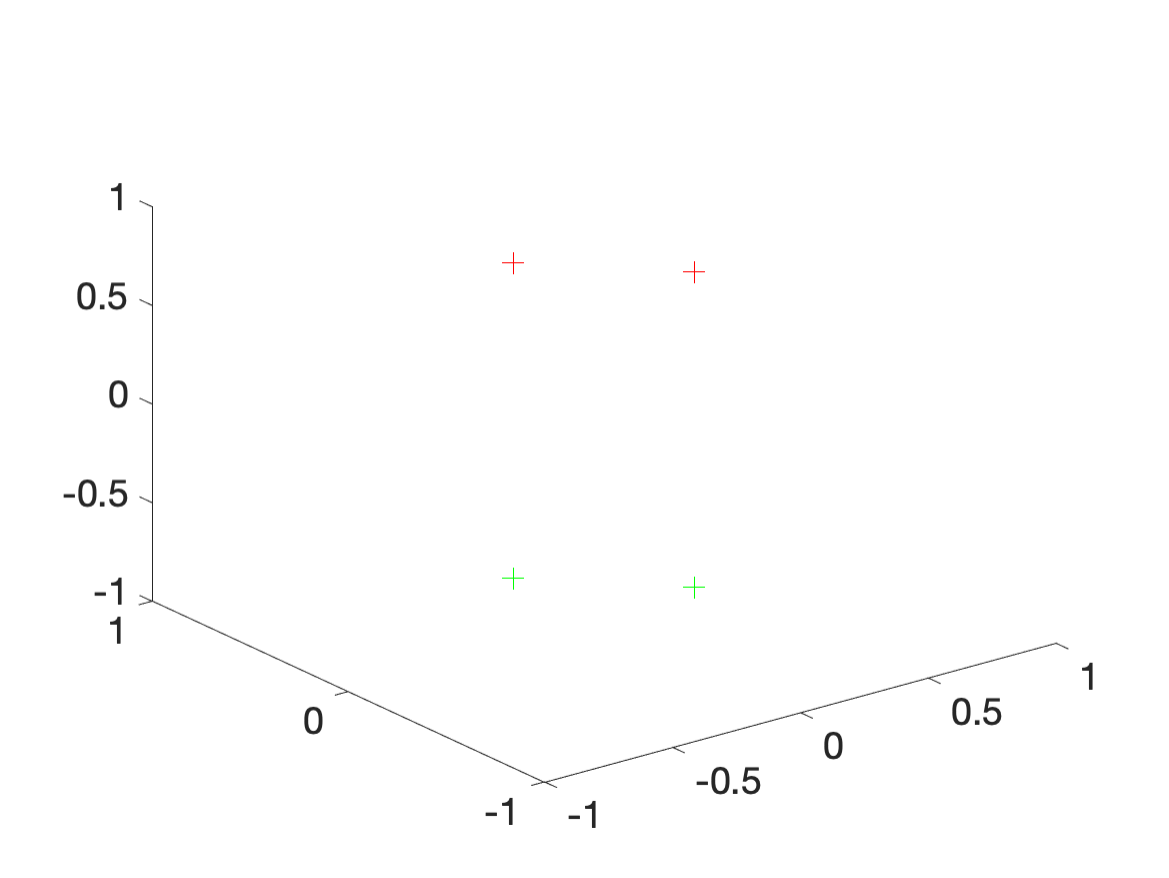}
  \includegraphics[scale=0.25]{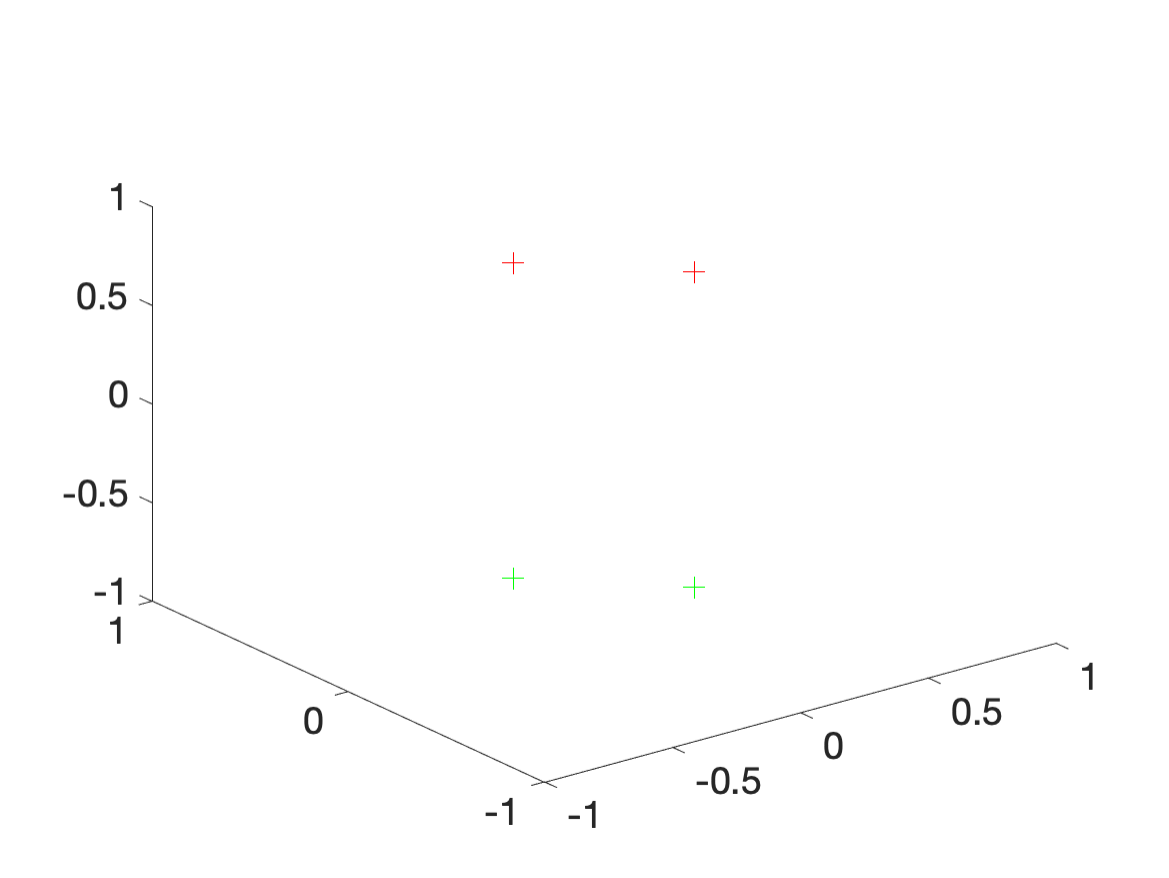} 
  \caption{Fourier inversion. $G(s,x) = \exp(\pi i s \cdot x)$. $X=[-1,1]^3$. $\{s_j\}$ are $J=8192$ randomly chosen points in $[-4,4]^3$. Columns: $\sigma$ equals to $10^{-3}$, $10^{-4}$, and $10^{-5}$, respectively. Rows: the easy case (top) and the hard case (bottom).}
  \label{fig:ex3F}
  \end{figure}

  Figure \ref{fig:ex3F} summarizes the experimental results. The three columns correspond to $\sigma$ equal to $10^{-3}$, $10^{-4}$, and $10^{-5}$. Two tests are performed. In the first one (top row), $\nx=4$ and $\{x_k\}$ are well-separated from each other.  In the second one (bottom row), the spikes again form two nearby pairs. The reconstructions are accurate for all $\sigma$ values.

\end{example}

\section{Discussions}\label{sec:disc}

This note extends the eigenmatrix approach to multidimensional unstructured sparse recovery problems. As a data-driven approach, it assumes no structure on the samples and offers a rather unified framework for such recovery problems. There are several directions for future work.
\begin{itemize}
\item A better understanding of the relationship between the size $n_a$ of the Chebyshev grid and the distribution of $\{s_j\}$.
\item Also how does the accuracy of the eigenmatrix $M$ depend on $\{s_j\}$ and the singular value threshold?
\item Providing the error estimates for the recovery problems mentioned in Section \ref{sec:intro}.
\item The recovery algorithm presented above follows the ESPRIT algorithm and \cite{andersson2018esprit}. An immediate question is whether other methods (such as the MUSIC and the matrix pencil method) can be extended and combined with the eigenmatrix approach.
\end{itemize}

\bibliographystyle{abbrv}

\bibliography{ref}

\end{document}